\documentclass{amsart}
\usepackage{amssymb}
\usepackage{amsfonts}
\usepackage{graphicx}
\usepackage{epstopdf}

\newcommand{\beq}{\begin{equation}}
\newcommand{\eeq}{\end{equation}}
\newtheorem{prop}{Proposition}
\newtheorem{thm}{Theorem}

\begin{document}

\title{A universal formula for the volume of compact Lie groups}
 
    \author{R. L. Mkrtchyan}
       \address{Yerevan Physics Institute, 2 Alikhanian Brothers St., Yerevan, 0036, Armenia}
           \email{mrl55@list.ru}
           
        \author{A. P. Veselov}
       \address{Department of Mathematics,
        Loughborough University, Loughborough,
        Leicestershire, LE11 3TU, UK and
       Moscow State University, Russia}
       \email{A.P.Veselov@lboro.ac.uk}

\maketitle




{\small  {\bf Abstract.} We provide a closed formula for the volume of a simple compact Lie group in terms of the universal Vogel parameters.  For the unitary  groups $SU_n$ this reduces to the integral representation of the classical Barnes G-function.}


\section{Introduction}

Let $G$ be a connected and simply-connected compact Lie group with simple Lie algebra $\mathfrak g.$ Supply $G$ with the Riemannian metric defined by the canonical Cartan-Killing form. What is the volume of $G$ ? 

Although some calculations for the classical groups were done already by Hurwitz and Weyl
the first general effective formula for the group volume seems to be found only in 1980 by Macdonald \cite{Ma}.
Macdonald's formula says that
\beq
\label{mac1}
Vol(G)=Vol(\mathfrak g/\mathfrak g_{\mathbb Z})\prod_{i=1}^r \frac{2\pi^{m_i+1}}{m_i!}.
\eeq
where $r$ is rank and $m_i,\, i=1,\dots, r$ are the exponents of the corresponding simple Lie algebra $\mathfrak g$,
$\mathfrak g_{\mathbb Z}$ is the Chevalley lattice (see \cite{Bour, Ma}).
An elegant topological proof of this formula was found by Hashimoto \cite{Hash}.
Different expressions for the group volume were proposed by Marinov \cite{Mar} and Kac and Peterson \cite{KP}.

The main result of this paper is the proof of the following formula for the group volume, 
which came naturally from our study of the universal Chern-Simons theory \cite{MV,Mkr}
(see formula (32) in \cite{Mkr}):
\beq
\label{main2}
Vol(G)=(2\sqrt 2\pi)^{\dim G} e^{-\Phi(\alpha,\beta,\gamma)},
\eeq
where 
\beq
\label{int}
\Phi(\alpha,\beta,\gamma)=\int_0^{\infty} \frac{F(x;\alpha,\beta,\gamma) dx}{x(e^x-1)},
\eeq
\beq
\label{gene}
F(x;\alpha,\beta,\gamma)=\frac{\sinh(x\frac{\alpha-2t}{4t})}{\sinh(\frac{x\alpha}{4t})}\frac{\sinh(x\frac{\beta-2t}{4t})}{\sinh(x\frac{\beta}{4t})}\frac{\sinh(x\frac{\gamma-2t}{4t})}{\sinh(x\frac{\gamma}{4t})}-\frac{(\alpha-2t)(\beta-2t)(\gamma-2t)}{\alpha\beta\gamma}.
\eeq
Here $t=\alpha+\beta+\gamma$, where $\alpha, \beta, \gamma$ are the universal parameters of the simple Lie algebra $\mathfrak g$ introduced by Vogel \cite{V}
and $\dim \, G$ is the dimension of $G$, which can be given by Vogel's formula
\beq
\label{dimG}
\dim G= \frac{(\alpha-2t)(\beta-2t)(\gamma-2t)}{\alpha\beta\gamma}.
\eeq
The parameters $\alpha, \beta, \gamma$ depend on the choice of an invariant bilinear form on $\mathfrak g$ and thus are defined up to common multiple and up to permutation, see \cite{V}. Their values are given in Table 1, where we have chosen the normalisation with $\alpha=-2$ corresponding to the so-called minimal bilinear form, when the square of the length of the maximal root equals 2 and $t=h^\vee$ is the dual Coxeter number. Note that the orthogonal groups are not simply connected with the fundamental group
$\pi_1(SO_n) = \mathbb Z_2$ for $n>2$, so they are represented by their double covers $Spin_n$ (having the volume twice bigger).

\begin{table}[h]  
\caption{Vogel's parameters for compact Lie groups}     
\begin{tabular}{|c|c|c|c|c|}
\hline
Lie group  & $\alpha$ & $\beta$ & $\gamma$  & $t=h^\vee$\\   
\hline    
$SU_n$     & $-2$ & 2 & $n $ & $n$\\
$Spin_{n}$    & $-2$ & 4& $n-4 $ & $n-2$\\
$Sp_{2n}$    & $-2$ & 1 & $n+2 $ & $n+1$\\
$G_{2}  $    & $-2$ & $10/3 $& $8/3$ & $4$ \\
$F_{4}  $    & $-2$ & $ 5$& $ 6$ & $9$\\
$E_{6}  $    & $-2$ & $ 6$& $ 8$ & $12$\\
$E_{7}  $    & $-2$ & $ 8$& $ 12$ & $18$ \\
$E_{8}  $    & $-2$ & $ 12$& $20$ & $30$\\
\hline  
\end{tabular}
\end{table}

The integral (\ref{int}) converges and defines an analytic function in the complement $ \mathbb CP^2 \setminus \Delta$ to the set 
$$\Delta = \{ (\alpha: \beta: \gamma) \in \mathbb CP^2: Re(\frac{\alpha}{t}) \geq 0,\, Re(\frac{\beta}{t}) \geq 0,\, Re(\frac{\gamma}{t}) \geq 0 \}.$$
We show that on the $SU_n$ line the corresponding function $\Psi(z)=\Phi(-2, 2, z)$ can be given by
\beq
\label{barnes}
\Psi(z)=\ln \,G(z+1) -\frac{1}{2} z^2 \ln\, z +\frac{1}{2}(z^2-z)\ln\, 2\pi,
\eeq
where $G$ is the classical Barnes' $G$-function \cite{B1}.
The appearance of $G$-function is not too surprising since it is a "natural" analytic continuation of the product of factorials
$$G(n+1)=1!2!3!\dots(n-1)!, \,\,\, n \in \mathbb Z_+,$$
which is an essential part of Macdonald's volume formula for $SU_n$:
\beq
\label{sun}
Vol(SU_n)=\frac{2^{\frac{n^2-1}{2}} n^{\frac{n^2}{2}} (2\pi)^{\frac{n^2+n-2}{2}}}{1!2!3!\dots(n-1)!}.
\eeq
However since such an analytic continuation is not unique, the fact that the analytic continuation 
provided by our formula agrees with the one proposed by Barnes, seems to be remarkable 
and shows a close relation of the volume function with the theory of multiple 
Barnes' gamma and zeta functions \cite{B2,B3}.

\section{Proof of the main formula}

We start with the following result expressing the group volume in terms of the root systems (see Kac-Peterson \cite{KP}, Prop. 4.32).
Let $\mathfrak g_{\mathbb C}$ be the complexification of the Lie algebra $\mathfrak g$, $\mathfrak h \subset \mathfrak g_{\mathbb C}$ be its 
Cartan subalgebra, $R_+$ be a positive part of the root system $R \in \mathfrak h^*$,
$$\rho= \frac{1}{2} \sum _{\mu \in R_+} \mu$$
be the Weyl vector and $< , >$ denote the Cartan-Killing form on $\mathfrak h^*$ 
\beq
\label{CK}
<X,Y>=tr \, ad_X ad_Y,\, X,Y \in \mathfrak g.
\eeq
It is known to be negative definite on $\mathfrak g$, so taken with sign minus determines the Riemannian metric on $G,$
which we call Cartan-Killing metric.

\begin{prop} (Kac and Peterson \cite{KP}).
The volume of a simple, connected and simply connected compact Lie group $G$ with Cartan-Killing metric can be given by
\beq
\label{KP}
Vol(G)=(2\sqrt {2} \pi)^{\dim\, G} \prod_{\mu \in R_+} \frac{\sin 2\pi <\rho, \mu>}{2\pi <\rho, \mu>}.
\eeq
 \end{prop}

Consider
$$\Phi= -\ln \prod_{\mu \in R_+} \frac{\sin 2\pi <\rho, \mu>}{2\pi <\rho, \mu>}$$
and rewrite it using the Euler formula
\begin{eqnarray}
\frac{\sin \pi x}{\pi x}=\frac{1}{\Gamma(1-x) \Gamma(1+x)},
\end{eqnarray}
where $\Gamma$ is the classical gamma function \cite{WW}:
$$
\Phi=\sum_{\mu \in R_+}\ln \Gamma(1-2<\rho,\mu>)+\ln(\Gamma(1+2<\rho,\mu>).
$$
Now we use Malmsten's integral representation of the logarithm of the gamma function, see \cite{WW}, Ch. 12.31:
$$
\ln\Gamma(1+z)=\int_{0}^{\infty}dx \frac{e^{-zx}+z(1-e^{-x})-1}{x(e^x-1)},
$$
implying that
$$
\ln\Gamma(1+z) + \ln\Gamma(1-z)=\int_{0}^{\infty}dx \frac{e^{zx}+e^{-zx}-2}{x(e^x-1)}.
$$
Thus we have
$$
\Phi=\sum_{\mu \in R_+}\int_{0}^{\infty}dx \frac{e^{2<\rho,\mu>x}+e^{-2<\rho,\mu>x}-2}{x(e^x-1)}=\int_{0}^{\infty}dx \frac{\sum_{\mu \in R}(e^{2<\rho,\mu>x}-1)}{x(e^x-1)}.
$$
Now we use the key relation 
$$
\sum_{\mu \in R}(e^{2<\rho,\mu>x}-1)=\frac{\sinh(x\frac{\alpha-2t}{4t})}{\sinh(\frac{x\alpha}{4t})}\frac{\sinh(x\frac{\beta-2t}{4t})}{\sinh(x\frac{\beta}{4t})}\frac{\sinh(x\frac{\gamma-2t}{4t})}{\sinh(x\frac{\gamma}{4t})}-\frac{(\alpha-2t)(\beta-2t)(\gamma-2t)}{\alpha\beta\gamma},
$$
see \cite{MV}, formula (2.9). This means that 
$$
\Phi=\int_0^{\infty} \frac{F(x;\alpha,\beta,\gamma) dx}{x(e^x-1)}
$$
with $F$ defined by (\ref{gene}). Thus we have proved

\begin{thm}
The volume of a simple, connected and simply connected compact Lie group $G$ 
with the canonical Cartan-Killing metric can be given in terms of Vogel's parameters by formula (\ref{main2}).
\end{thm}

One can easily check that the integral (\ref{int}) converges for all complex values of parameters $\alpha, \beta, \gamma$ unless
$Re(\frac{\alpha}{t}) \geq 0,\, Re(\frac{\beta}{t}) \geq 0,\, Re(\frac{\gamma}{t}) \geq 0$ simultaneously. In particular, this is the case when all parameters are real and not of the same sign, which is is true for the parameters corresponding to simple Lie algebras.

In the next section we show that on the unitary line $\alpha+\beta=0$, containing all the unitary groups, our integral representation essentially reduces to the one for the classical Barnes $G$-function \cite{B1}. The relations with the rich theory of the multiple Barnes gamma and zeta functions \cite{B2,B3} should help to study the analytic properties of the functions $\Phi(\alpha,\beta,\gamma)$ and $e^{\Phi(\alpha,\beta,\gamma)}$, which will be done elsewhere.

\section{Volume of the unitary groups and the Barnes $G$-function}

Vogel's parameters of $G=SU_n$ are $\alpha=-2, \beta=2, \gamma=n.$ Replacing here $n \in \mathbb N$ by $z$ (real or complex) we get the line on the Vogel plane defined by $\alpha+\beta=0.$ A simple calculation shows that the restriction of the function (\ref{gene}) to this line is
\beq
\label{su}
F(x;-2,2,z)=\frac{\cosh x-1}{2\sinh^2 \frac{x}{2z}}-z^2.
\eeq
Since
$$\frac{\cosh x-1}{e^x-1}=\frac{1-e^{-x}}{2}$$
we have 
$$
\Psi(z):=\Phi(-2,2,z)=\int_{0}^{\infty}\left(\frac{1-e^{-x}}{4\sinh^2 \frac{x}{2z}}-\frac{z^2}{e^x-1}\right)\frac{dx}{x}=\int_{0}^{\infty}\left(\frac{1-e^{-zy}}{4\sinh^2 \frac{y}{2}}-\frac{z^2}{e^{zy}-1}\right)\frac{dy}{y}
$$
or, after elementary transformation,
$$
\Psi(z)=-\int_{0}^{\infty}\left(\frac{e^{-(z+1)y}}{(1-e^{-y})^2}+\frac{z^2}{e^{zy}-1}-\frac{e^{-y}}{(1-e^{-y})^2}\right)\frac{dy}{y}.
$$
Compare this with the integral representation of the logarithm of the Barnes $G$-function:
\beq
\label{G}
\ln G(z+1)=\frac{z+1}{2}\ln 2\pi-\int_{0}^{\infty}\left(\frac{e^{-(z+1)y}}{(1-e^{-y})^2}-\frac{1}{y^2}+\frac{z}{y}-\frac{e^{-y}}{2}(z^2-\frac{1}{6})\right)\frac{dy}{y},
\eeq
see Barnes \cite{B2}, pages 328 and 368.
Consider the difference
$$
\Psi(z)-\ln G(z+1)= -\frac{z+1}{2}\ln 2\pi+I_1 + I_2,
$$
where 
$$I_1=\int_{0}^{\infty}\left(\frac{e^{-y}}{(1-e^{-y})^2}-\frac{1}{y^2}+\frac{1}{12}e^{-y}\right)\frac{dy}{y}, \, 
I_2=\int_{0}^{\infty}\left(\frac{z}{y}-\frac{z^2}{e^{zy}-1}-\frac{e^{-y}}{2}z^2\right)\frac{dy}{y}.$$
The first integral is a constant:
$I_1=\frac{1}{2}\ln 2\pi,$ as it follows from (\ref{G}) since $G(1)=1$. 
Change of variable $y \rightarrow t=zy$ in the second integral gives
$$
I_2=\int_{0}^{\infty}\left(\frac{z^2}{t}-\frac{z^2}{e^{t}-1}-\frac{e^{-t/z}}{2}z^2\right)\frac{dt}{t}=z^2\int_{0}^{\infty}\left(\frac{1}{t}-\frac{1}{e^{t}-1}-\frac{e^{-t}}{2}\right)\frac{dt}{t}$$
$$+\frac{z^2}{2}\int_{0}^{\infty}(e^{-t}-e^{-t/z})\frac{dt}{t} = Cz^2-\frac{1}{2} z^2 \ln z
$$
where we have used Frullani's integral
$$
\int_{0}^{\infty}(e^{-pt}-e^{qt})\frac{dt}{t}=\ln \frac{q}{p},
$$
see e.g. Gradshtein-Ryzhik \cite{GR}, formula 3.434.2. The constant
$$C=\int_{0}^{\infty}\left(\frac{1}{t}-\frac{1}{e^{t}-1}-\frac{e^{-t}}{2}\right)\frac{dt}{t}=\frac{1}{2}\ln 2\pi,$$
since $G(2)=1$ in (\ref{G}). Thus we have

\begin{thm}
On the line $\alpha+\beta=0$ on the Vogel plane 
containing the unitary groups $SU_n$ the function $\Phi(-2,2,z)$ defined by (\ref{int}) gives an analytic continuation of $\Phi(SU_n)$, 
which is related to the classical Barnes $G$-function by 
\beq
\label{barnes2}
\Phi(-2,2,z)=\ln \,G(z+1) -\frac{1}{2} z^2 \ln\, z +\frac{1}{2}(z^2-z) \ln\,2\pi.
\eeq
\end{thm}

Let's compare this with the explicit formula (\ref{sun}) for the volume of $SU_n$. From (\ref{sun}) we have since $\dim\, SU_n=n^2-1$ that
$$
\frac{Vol(SU_n)}{(2\sqrt {2} \pi)^{\dim\, SU_n}}=\frac{n^{\frac{n^2}{2}} (2\pi)^{\frac{n-n^2}{2}}}{1!2!3!\dots(n-1)!}.
$$
or,
\beq
\label{sun1}
\ln \frac{Vol(SU_n)}{(2\sqrt {2} \pi)^{\dim\, SU_n}}=\frac{n^2}{2} \ln\, n  -\frac{1}{2}(n^2-n)\ln\,2\pi -\ln(1!2!3!\dots(n-1)!).
\eeq
Comparing this with (\ref{barnes2}) we see that that our integral representation gives the analytic continuation 
of the product of the factorials $1!2!3!\dots(n-1)!$ coinciding with the one proposed by Barnes \cite{B1}.
At the level of the asymptotic expansions this agreement was shown in \cite{Mkr}.

We should say that the Barnes analytic continuation is widely used and plays an important role in modern theoretical physics, 
in particular, in topological string and Chern-Simons theory, see e.g. Periwal \cite{Pe} and Ooguri and Vafa \cite{OV}. 
Our universal volume formula may have some applications here as well.


\section{Acknowledgments.}

This work is partially supported by the grant of Science Committee of Ministry of Science and Education of the Republic of Armenia, 
Volkswagen Foundation, and the EPSRC (grant EP/J00488X/1).

\end{document}